\documentclass [10pt,a4paper,draft] {article}
\usepackage [applemac] {inputenc}
\usepackage[french] {babel}
\usepackage {amsfonts}
\usepackage {amsmath}
\usepackage {amssymb}

\begin{document}

\begin{center}
\begin{LARGE}
Les invariants polyn\^omes de la repr\'esentation 

\smallskip
coadjointe de groupes inhomog\`enes

\end{LARGE}
\bigskip

par : Mustapha RAIS  (Poitiers)
%\footnotemark \footnotetext[1]{Universit\'e de POITIERS -
%D\'epartement de Math\'ematiques - T\'el\'eport 2, Boulevard Marie et Pierre Curie - BP 30179 -
%86962 FUTUROSCOPE CHASSENEUIL Cedex} 
\end{center}

\bigskip
\vspace{10 mm}

\begin{large}
\noindent
\textbf{        Introduction}
\end{large}

\medskip
        R\'ecemment Akaki TIKARADZE m'a pos\'e diverses questions sur le centre $ZU(\mathfrak{p})$ de l'alg\`ebre
enveloppante $U(\mathfrak{p})$ d'alg\`ebres de Lie $\mathfrak{p} = \mathfrak{g} \underset{\rho}{\times} V$,
produits semi-directs d'une alg\`ebre de Lie $\mathfrak{g}$ par un espace vectoriel $V$, relativement \`a une
repr\'esentation lin\'eaire $\rho$ de $\mathfrak{g}$ dans $V$. A isomorphisme d'alg\`ebres pr\`es
(l'isomorphisme de Duflo), on peut se contenter d'examiner l'alg\`ebre $Y(\mathfrak{p})$ des invariants de
$\mathfrak{p}$ dans l'alg\`ebre sym\'etrique $S(\mathfrak{p})$ de $\mathfrak{p}$ (i.e. dans l'alg\`ebre des
fonctions polyno\-miales sur l'espace vectoriel $\mathfrak{p}^*$ des formes lin\'eaires sur $\mathfrak{p}$). Il se
trouve que les alg\`ebres de Lie $\mathfrak{p}$ auxquelles s'int\'eressait A. Tikaradze ont donn\'e lieu \`a
plusieurs travaux anciens et r\'ecents. Le lecteur trouvera ci-dessous dans plusieurs notes bibliographiques des
indications pr\'ecises sur certains de ces travaux, en particulier ceux de Hitoshi KANETA (\cite{Ka-1} et
\cite{Ka-2}) et de Dmitri PANYUSHEV (\cite{Pa}). On trouvera dans les bibliographies des articles de Kaneta et de
Panyushev une liste de travaux portant sur ces questions, que je n'ai pas consult\'es et que je m'excuse de ne pas
citer.

\smallskip
        Le but de ce texte, qui ne pr\'etend pas \`a une grande originalit\'e, est d'expliciter pour certaines alg\`ebres
d'invariants dans $S(\mathfrak{p})$, un syst\`eme de g\'en\'erateurs.

\bigskip
\vspace{5 mm}

\begin{large}
\noindent
\textbf{Notations :}
\end{large}

- Le corps de base est $\mathbb{C}$. L'espace vectoriel $V$ est $\mathbb{C}^n$, identifi\'e \`a l'espace $M_{n,1}$ des
matrices colonnes, son dual $V^*$ est identifi\'e \`a l'espace $M_{1,n}$ des matrices lignes.

\vskip 3mm
- Lorsque $x$ est une matrice carr\'ee de taille $n$, on note $p_1(x),\ldots , p_n(x)$ les coefficients du
polyn\^ome caract\'eristique \'ecrit sous la forme :
$$
        det(tI_n-x) = t^n - p_1(x)t^{n-1} - p_2(x)t^{n-2}- \cdots - p_n(x)
$$
et $B_0,B_1,\ldots ,B_{n-1}$ les gradients respectifs de $p_1,\ldots ,p_n$ calcul\'es au moyen de la forme trace :
$$
        tr(B_k(x)y) = \big(\frac{d}{dt}\big)_0 \ p_{k+1}(x+ty) \quad (y \in M_n(\mathbb{C})).
$$
Les expressions explicites des $B_k$ sont les suivantes :
\begin{eqnarray}
B_0(x) &= &I_n \nonumber\\
        B_k(x)                &= & x^k - p_1(x)x^{k-1}- \cdots - p_k(x)I_n \quad (1\leq k \leq n-1)  \nonumber\\
\nonumber
\end{eqnarray}

        Les cas trait\'es ici sont ceux des groupes ``inhomog\`enes'' :
$$
        SL(n) \times \mathbb{C}^n = ISL(n), \quad O(n) \times \mathbb{C}^n = IO(n),
$$
$$
        SO(n) \times \mathbb{C}^n = ISO(n)
$$
et celui du groupe : 

$$
        GL(n) \underset{\rho}{\times} (V \oplus V^*) \quad (\rho(g).(u,v^*) = (g.u, v^*g^{-1})
$$
qui est tr\`es peu diff\'erent d'une $\mathbb{Z}_2$-contraction au sens de Panyushev et qui est utilis\'e dans 
\cite{Ti}.

%\vfill\eject
\vskip 10mm
\begin{large}
\noindent
\textbf{1~Le cas $ISL(n)$}
\end{large}

%%%%%%%%%%%%%%%%%%%%%%%
%\vskip 10mm
%\section{Le cas $ISL(n)$}
        
\vskip 7mm
\noindent
\textbf{1.1. } Au pr\'ealable, on revient sur le groupe affine de $\mathbb{C}^n : A = G \underset{\rho}{\times}
\mathbb{C}^n$, avec $G = GL(n)$ et $\rho(g).u = g.u$, dont la repr\'esentation coadjointe a \'et\'e \'etudi\'ee en
d\'etail  dans \cite{Ra-1}. L'alg\`ebre de Lie $\mathfrak{a} = \mathfrak{g}\underset{\rho}{\times} V$,
(o\`u $\mathfrak{g} =
\mathfrak{g}\ell(n))$, est d'indice z\'ero et l'alg\`ebre $Y'(\mathfrak{a})$ engendr\'ee par les semi-invariants
de $A$ dans $S(\mathfrak{a})$ est engendr\'ee par un polyn\^ome $f$, dont on rappelle l'expression
ci-dessous.

\vskip 1mm
        Pour cela, on identifie $\mathfrak{a}^*$ \`a l'espace vectoriel $\mathfrak{g} \times V^*$ : un couple $(y,v^*)$
dans $\mathfrak{g} \times V^*$ d\'efinit sur $\mathfrak{a}$ la forme lin\'eaire $\ell$ telle que, pour tout
$(x,u)$ dans $\mathfrak{a}$ :
$$
        <\ell, (x,u)> = tr(yx) + v^*u
$$
(o\`u $v^*u$ est le produit matriciel, et $v^*u = tr(uv)$ si on pr\'ef\`ere).

\vskip 1mm
        L'alg\`ebre $Y'(\mathfrak{a})$ est $\mathbb{C}[f]$, o\`u $f$ est d\'efini comme un d\'eterminant de $n$ vecteurs
lignes :
$$
        f(y,v^*) = det(v^*B_{n-1}(y), v^*B_{n-2}(y),\ldots , v^*)
$$
pour tout $(y,v^*)$ dans $\mathfrak{a}^*$, et on a, lorsque $(g,u) \in A$ : 
$$
        f(Ad^*(g,u).(y,v^*)) = (det \ g)^{-1} f(y,v^*).
$$
\indent
Tout ceci est d\'emontr\'e dans \cite{Ra-1} (\textit{c.f.} en cas de besoin les num\'eros 2.15, 3.7 et 3.8).

\vskip 7mm
\noindent
\textbf{1.2.}~On s'int\'eresse maintenant aux covariants $\Phi : \mathfrak{a}^* \longrightarrow V^*$,
c'est-\`a-dire aux fonctions polynomiales $\Phi$ telles que :
$$
        \Phi(Ad^*(g,u).\ell) = \Phi(\ell)g^{-1} \qquad ((g,u) \in A,\ \  \ell \in \mathfrak{a}^*).
$$
On d\'efinit $n$ fonctions $\Phi_k$ en posant :
$$
        \Phi_k(y,v^*) = v^*B_k(y)\qquad (0 \leq k \leq n-1)
$$

\bigskip
\noindent
\textbf{Lemme} : (\cite{Ra-1}, 4.6) 
\vskip 1mm
\indent
        \textit{ L'espace vectoriel des covariants $\Phi : \mathfrak{a}^* \longrightarrow V^*$ est de dimension $n$ et
$(\Phi_0,\ldots , \Phi_{n-1})$ en est une base.}

\vskip 7mm
\noindent
\textbf{1.3.}~Soit $P = ISL(n)$ le sous-groupe de $A$ constitu\'e par les couples $(g,u)$ avec : $det\ g = 1$. Son
alg\`ebre de Lie $\mathfrak{p}$ s'identifie \`a la sous-alg\`ebre de Lie de $\mathfrak{a}$ constitu\'ee par les
couples $(y,v^*)$ tels que $tr(y) = 0$, de sorte que $\mathfrak{p} = [\mathfrak{a}, \mathfrak{a}]$ est
l'alg\`ebre d\'eriv\'ee de $\mathfrak{a}$. D'apr\`es \cite{Ra-1} ($\no 3.10$), l'alg\`ebre $Y(\mathfrak{p})$ est
$\mathbb{C}[\overline{f}]$, o\`u $\overline{f}$ est la restriction de $f$ \`a $\mathfrak{p}^*$, $\mathfrak{p}^*$
\'etant identifi\'e au sous-espace $sl(n) \times V^*$ de $\mathfrak{a}^*$.

\vskip 7mm
\noindent
\textbf{1.4.}~\underline{Note bibliographique} :
\vskip 1mm
        L'\'etude du groupe $ISL(n)$ appara\^it dans un article de Hitoshi Kaneta (\cite{Ka-1}). L'auteur introduit un
sous-espace vectoriel de dimension $n$ de $\mathfrak{p}^*$, not\'e $\mathfrak{h}$, qui est l'espace des
$(y,v^*)$ avec :
$$
        y = a_1\, E_{21} +\cdots + a_{n-1}\, E_{n,n-1}, \qquad v^* = be^*_n
$$
o\`u $(E_{ij})$ est la base canonique de $\mathfrak{g}\ell(n),\ e^*_n = (0,\ldots , 0,1)$, $(a_1,\ldots ,
a_{n-1}, b) \in \mathbb{C}^n$, d\'efinit une fonction polyn\^ome $t$ sur $\mathfrak{h}$ : 
$$
        t(a_1E_{21} +\cdots + a_{n-1}\, E_{n,n-1}, be^*_n) = (\prod_{1\leq k\leq n-1} \ a^k_k)b^n
$$
et \'enonce le :

\bigskip
\noindent
\textbf{Th\'eor\`eme} : (\cite{Ka-1}, Theorem 1) 

        \textit{L'application de restriction de $Y(\mathfrak{p})$ dans l'ensemble des polyn\^omes sur $\mathfrak{h}$
est un homomorphisme d'alg\`ebres injectif, dont l'image est $\mathbb{C}[t]$.}

\vskip 3mm
        On constate alors que $t$ n'est autre que la restriction de $\overline{f}$ \`a $\mathfrak{h}$. La contribution
de [Ra-1] consiste donc en l'expression explicite, en fonction des coordonn\'ees naturelles
$(y,v^*)$, du g\'en\'erateur de $Y(\mathfrak{p})$ dont la restriction \`a $\mathfrak{h}$ est $t$. On notera
que \cite{Ra-1} date de 1978, alors que \cite{Ka-1} date de 1984.

\vskip 10mm
\begin{large}
\noindent
\textbf{2~Le cas $GL(n)\underset{\rho}{\times} (V \oplus V^*)$}
\end{large}

\bigskip 
\noindent
\textbf{2.1.}~Ici aussi, on pose $GL(n) = G$ et $\mathfrak{g}\ell(n) = \mathfrak{g}$ ; la repr\'esentation $\rho$
est la somme directe de la repr\'esentation naturelle de $G$ dans $V$ et de sa repr\'esentation duale dans
$V^*$ :
$$
        \rho(g).(u,v^*) = (g.u,\, v^*g^{-1})
$$
\indent
On notera $B$ le groupe consid\'er\'e, $\mathfrak{b}$ son alg\`ebre de Lie, et $\mathfrak{b}^*$ le dual de
$\mathfrak{b}$. Comme espace vectoriel, on a : $\mathfrak{b} = \mathfrak{g} \times V \times V^*$, avec le
crochet qui d\'efinit le produit semi-direct $\mathfrak{g} \underset{\rho'}{\times} (V \oplus V^*)$ : 
$$
        [(x,0,0), (0,u,v^*)] = (0, x.u, -v^*x)
$$
$$
        \rho'(x).(u,v^*) = (x.u, -v^*x)
$$
\indent
Comme espace vectoriel, on a : $\mathfrak{b}^* = \mathfrak{g} \times V^* \times V$ ; un \'el\'ement 
$(y,w^*,\xi)$ de $\mathfrak{g} \times V^* \times V$ d\'efinit la forme lin\'eaire $\ell$ sur $\mathfrak{b}$,
telle que, pour tout $(x,u,v^*)$ dans $\mathfrak{b}$ : 
$$
        <\ell, (x,u,v^*)> \, = tr(yx) + w^*u + v^*\xi
$$
\indent
        On identifie $G$ au sous-groupe $G \times (0,0)$. Lorsque $g \in G$, on a alors : 
$$
        Ad_B(g).(x,u,v^*) = (gxg^{-1}, g.u, v^*g^{-1})
$$
$$
        Ad^*_B(g).(y,w^*,\xi) = (gyg^{-1}, w^*.g^{-1},g.\xi)
$$

\vskip 7mm
\noindent
\textbf{2.2.}~L'application $(g,u) \longmapsto (g,u,0)$ identifie le groupe $A$ qui a fait l'objet du \no 1
ci-dessus, \`a un sous-groupe de $B$, ce qui induit une identification de l'alg\`ebre de Lie $\mathfrak{a}$ \`a
la sous-alg\`ebre de $\mathfrak{b}$ constitu\'ee par les triplets $(x,u,0)$, et de m\^eme le dual
$\mathfrak{a}^*$ au sous-espace de $\mathfrak{b}^*$ constitu\'e par les $(y,w^*,0)$.

        Ceci \'etant, on constate que $B$ est le produit semi-direct de $A$ par $V^*$, relativement \`a la
repr\'esentation lin\'eaire $\delta$ de $A$ dans $V^*$ : 
$$
        \delta(g,u).v^* = v^*g^{-1}
$$
et que $\mathfrak{b}$ est le produit semi-direct de l'alg\`ebre de Lie $\mathfrak{a}$ par $V^*$ relativement
\`a la repr\'esentation lin\'eaire $\delta'$ de $\mathfrak{a}$ dans $V^*$ avec : 
$$
        \delta'(x,u).v^* = -v^*x.
$$
\indent
        L'indice de la repr\'esentation $\delta'$ est z\'ero et l'application de la formule de l'indice (\cite{Ra-2}) montre
que l'indice de l'alg\`ebre $\mathfrak{b}$ est $n$.

        On utilisera les r\'esultats de \cite{Ra-1}, dont certains ont \'et\'e rappel\'es ci-dessus, pour d\'eterminer de
fa\c con explicite un syst\`eme de g\'en\'erateurs de $Y(\mathfrak{p})$.

\vskip 7mm
\noindent
\textbf{2.3.} Dans $\mathfrak{a}^*$, l'orbite coadjointe ouverte est l'ensemble $\Omega$ des $(y,w^*)$ tels que
$f(y,w^*) \not= 0$. Parmi les formes r\'eguli\`eres sur $\mathfrak{a}$, on a privil\'egi\'e la forme $\ell_0 =
(J,e^*_n)$, o\`u $e^*_n = (0,\ldots ,0,1)$ et $J$ est la matrice nilpotente (principale) d\'efinie par :
$$
        e^*_1 J = 0,\quad e^*_2 J = e^*_1,\ldots , e^*_nJ = e^*_{n-1}
$$
($(e^*_1,\ldots ,e^*_n)$ est la base naturelle de $\mathbb{C}^n = M_{1,n}$). Le stabilisateur de $\ell_0$ dans $A$
est trivial (\cite{Ra-1}, 3.2) et l'application bijective orbitale $T : A \longrightarrow \Omega\ (T(a) =
Ad^*_A(a).\ell_0)$ admet une application r\'eciproque $T^{-1}$ qui a \'et\'e d\'etermin\'ee explicitement
(\cite{Ra-1}, 3.8).

        On note $\widetilde{\Omega}$ l'ouvert $\Omega \times V$ de $\mathfrak{b}^*$. Soit $\ell = (y, w^*,\xi)$
dans $\widetilde{\Omega}$. L'unique \'el\'ement $(g,u)$ de $A$ tel que $Ad^*_A(g,u).\ell = \ell_0$ est tel que :
$$
        g= L(w^*B_{n-1}(y),\ w^ * B_{n-2}(y),\ldots , w^*)
$$
o\`u on a mis en \'evidence les lignes de $g$. Pr\'ecis\'ement : 
$$
        e^*_1g = w^*B_{n-1}(y),\ e^*_2g = w^*B_{n-2}(y),\ldots , e^*_ng = w^*
$$
\indent
Soit $F : \mathfrak{b}^* \longrightarrow \mathbb{C}$ une fonction invariante sous l'action coadjointe de $B$ dans
$\mathfrak{b}^*$. Lorsque $(y,w^*,\xi) = \ell$ appartient \`a $\widetilde{\Omega}$, on a donc :
$$
        F(\ell) = F(Ad^*_B(g,u).\ell) = F(J,e^*_n, g.\xi)
$$
o\`u $g$ est la matrice explicit\'ee ci-dessus, de sorte que $g.\xi$ est le vecteur colonne dont les coordonn\'ees
dans la base canonique $(e_1,e_2,\ldots , e_n)$ sont : $w^*B_{n-1}(y)\xi,\ w^*B_{n-2}(y)\xi,\ldots , w^*\xi$.
Pr\'ecis\'ement : 
$$
        g.\xi = \sum_k (w^*B_{n-k}(y).\xi)e_k
$$

\vskip 7mm
\noindent
\textbf{2.4.} \textbf{Proposition} :  

\smallskip
\begin{enumerate}
\item Les fonctions $F_k : \mathfrak{b}^* \longrightarrow \mathbb{C}$ d\'efinies par 
$$
        F_k(y,w^*,\xi) = w^*B_k(y)\xi \qquad (0\leq k \leq n-1)
$$
sont des fonctions (polynomiales) $Ad^*(B)$-invariantes sur $\mathfrak{b}^*$, alg\'ebriquement
ind\'epen\-dantes.

\medskip
\item $Y(\mathfrak{b}) = \mathbb{C}[F_0,\ldots , F_{n-1}]$ est l'alg\`ebre de polyn\^omes engendr\'ee par :
$F_0,\ldots , F_{n-1}$.
\end{enumerate}

\vfill\eject
%\medskip
\noindent
\textsc{D\'emonstration} : 
\begin{enumerate}
\item - Les fonctions $F_k$ sont li\'ees aux covariants $\Phi_k$ d\'efinis dans 1.2 de la mani\`ere suivante :
$$
        F_k(y,w^*,\xi) = \Phi_k(y,w^*)\xi
$$
\indent
Lorsque $(g,u)$ est un \'el\'ement de $A$, on a :

\begin{eqnarray}
F_k(Ad^*_B(g,u).(y,w^*,\xi)) &=                &F_k((Ad^*_A(g,u).(y,w^*),g\xi) \nonumber\\
        &=                &\Phi_k(Ad^*_A(g,u).(y,w^*)).g\xi\nonumber\\
&=        &\Phi_k(y,w^*)g^ {-1}g\xi\nonumber\\
&= & F_k(y,w^*,\xi)\nonumber
\end{eqnarray}

Ainsi les $F_k$ sont $Ad^*_B(A)$-invariantes. Pour montrer qu'elles sont $Ad^*(B)$-invariantes, il reste \`a
prouver que :
$$
        F_k(y + \xi v^*,w^*,\xi) = F_k(y,w^*,\xi)
$$
pour tout $v^*$ dans $V^*$, ou encore que :
$$
        B_k(y + \xi v^*)\xi = B_k(y)\xi \qquad (v^* \in V^*)
$$
On constate alors que ces \'egalit\'es r\'esultent imm\'ediatement des \'egalit\'es : $v^*B_k(y+uv^*) =
v^*B_k(y)$ d\'emontr\'ees dans (\cite{Ra-1}, 4.5).

\medskip
- Par ailleurs : $F_k(J,e^*_n,\xi) = e^*_nB_k(J)\xi = e^*_nJ^k\xi = \xi_{n-k}$ (avec : $\sum_k \xi_k\, e_k = \xi)$.
Ceci montre que les $F_k$ sont alg\'ebriquement ind\'ependantes.

\medskip
\item Soit $F : \mathfrak{b}^* \longrightarrow \mathbb{C}$ une fonction polynomiale $Ad^*(B)$-invariante et soit
$\overline{F} : V \longrightarrow \mathbb{C}$ d\'efinie par :
$$
        \overline{F}(\xi) = F(J,e^*_n,\xi) \qquad (\xi \in V)
$$
La fonction $\overline{F}$ est polynomiale et :
$$
        F(y,w^*,\xi) = \overline{F} (\sum_k\ F_k(y,w^*,\xi)e_{n-k})
$$
pour tout $(y,w^*,\xi)$ dans $\widetilde\Omega$, donc pour tout $(y,w^*,\xi)$ dans $\mathfrak{b}^*$.
\end{enumerate}

\vskip 7mm
\noindent
\textbf{2.5.}~Soit $\ell$ dans $\widetilde\Omega$. Le stabilisateur de $\ell$ dans $A$, \textit{i.e.} $B^\ell \ \cap
\ A$, est trivial. L'orbite sous $A$ de $\ell$, \textit{i.e.} $Ad^*_B(A).\ell$, est de codimension $n$ dans
$\mathfrak{b}^*$. Donc $\ell$ est une forme r\'eguli\`ere sur $\mathfrak{b}$. Il est s?r qu'il existe d'autres
formes r\'eguli\`eres sur $\mathfrak{b}$, par exemple les $\ell = (y,w^*,v)$ telles que les vecteurs :
$$
        v,B_1(y)v,\ldots , B_{n-1}(y)v
$$
soient lin\'eairement ind\'ependants. La d\'etermination de toutes les formes r\'eguli\`eres reste \`a faire.

\vskip 7mm
\noindent
\textbf{2.6.}~Soit $\pi : \widetilde\Omega \longrightarrow V$ d\'efinie par :
$$
        \pi(\ell) = \sum_k \ F_k(\ell)e_{n-k}
$$
Il est facile de voir que $\pi(\ell_1) = \pi(\ell_2)$ si et seulement si $\ell_1$ et $\ell_2$ sont
$Ad^*_B(A)$-conjugu\'es. Autrement dit : $\widetilde\Omega$ est un $A$-fibr\'e principal. De plus, lorsque
$\ell \in \widetilde\Omega$, on a : $(Ad^*_B(B).\ell)\ \cap \ \widetilde\Omega = Ad^*_B(A).\ell$ (l'orbite sous
$A$ de $\ell$ est un ouvert de l'orbite (coadjointe) sous $B$ de $\ell$).

\vskip 7mm
\noindent
\textbf{2.7.}~L'alg\`ebre $S(\mathfrak{b})^{\mathfrak{g}}$ des invariants de $\mathfrak{g} =
\mathfrak{g}\ell(n)$ dans $S(\mathfrak{b})$ est connue (\textit{c.f.} \cite{Ra-1}, \no 5.6). C'est l'alg\`ebre
engendr\'ee par les $2n$ fonctions (de $y,w^*,\xi)$ :
$$
        p_1(y),\ldots , p_n(y),\ w^*\xi, \ w^*y\xi,\ldots , w^*y^{n-1}\xi 
$$
ou, ce qui revient au m\^eme, par les fonctions :
$$
        p_1,p_2,\ldots , p_n, F_0, F_1,\ldots , F_{n-1}.
$$
On voit alors que le centralisateur (ou le commutant) de $\mathfrak{g}$ dans l'alg\`ebre enveloppante
$\cal{U}(\mathfrak{b})$ de $\mathfrak{b}$, est engendr\'e comme alg\`ebre par la r\'eunion du centre de
$\cal{U}(\mathfrak{g})$ et du centre de $\cal{U}(\mathfrak{b})$. Il est donc commutatif (voir \cite{Ra-3} pour
des r\'esultats de ce type).

\vskip 7mm
\noindent
\textbf{2.8.}~Il est possible de trouver une autre expression des $F_k$. Pour cela, on s'int\'eresse aux matrices
carr\'ees $X$ de taille $(n+1)$, qu'on \'ecrit sous la forme : 
$X = \begin{pmatrix} 
y        &v\\
w^*        &a\\
\end{pmatrix}
$, avec $y$ dans $\mathfrak{g}\ell(n)$, $v$ dans $V$, $w^*$ dans $V^*$ et $a$ dans $\mathbb{C}$. On peut alors
(\cite{Ra-3}, \no 1.3, formules (1;6)) calculer le polyn\^ome caract\'eristique de $X$, et on trouve :

\begin{eqnarray}
p_1(X) &=                &p_1(y)+a \nonumber\\
        p_{k+2}(X) &=                &p_{k+2}(y) - ap_{k+1}(y) + w^*B_k(y)v \quad (0 \leq k \leq n-2)\nonumber\\
p_{n+1}(X) &=        &-ap_n(y) + w^*B_{n-1}(y)v\nonumber
\end{eqnarray}

On a donc : 

\begin{eqnarray}
F_k(y,w^*,\xi) &=                &p_{k+2}\Bigg(\begin{pmatrix} 
y        &v\\
w^*        &0\\
\end{pmatrix}\Bigg) - p_{k+2}(y) \quad (0 \leq k \leq n-2) \nonumber\\
F_{n-1}(y,w^*,\xi) &=        &p_{n+1}\Bigg(\begin{pmatrix} 
y        &v\\
w^*        &0\\
\end{pmatrix}\Bigg) \nonumber
\end{eqnarray}

\vskip 7mm
\noindent
\textbf{2.9.}~\underline{Note bibliographique} :  

\smallskip
        La d\'ecomposition : $\begin{pmatrix} 
x        &u\\
v^*        &a\\
\end{pmatrix} = \begin{pmatrix} 
x        &0\\
0        &a\\
\end{pmatrix} + \begin{pmatrix} 
0        &u\\
v^*        &0\\
\end{pmatrix}$
fait de $\mathfrak{g}\ell(n+1)$ une alg\`ebre de Lie sym\'etrique : 
$$
        \mathfrak{g}\ell(n+1) = \mathfrak{g}_0 \oplus \mathfrak{g}_1,\ \mathfrak{g}_0 = \mathfrak{g}\ell(n) \times
\mathbb{C} = \Bigg\{ \begin{pmatrix} 
x        &0\\
0        &a\\
\end{pmatrix}\Bigg\}, \ \mathfrak{g}_1 = \Bigg\{ \begin{pmatrix} 
0        &u\\
v^*        &0\\
\end{pmatrix}\Bigg\}
$$
\`a laquelle est associ\'ee l'alg\`ebre de Lie $\mathfrak{k} = \mathfrak{g}_0 \times \mathfrak{g}_1$, produit
semi-direct de l'alg\`ebre de Lie $\mathfrak{g}_0$ par le $\mathfrak{g}_0$-module $\mathfrak{g}_1$. On se
reportera \`a \cite{Pa} pour une \'etude compl\`ete des alg\`ebres de Lie de ce type (les
$\mathbb{Z}_2$-contractions des alg\`ebres de Lie r\'eductives), l'alg\`ebre de Lie $\mathfrak{k}$ d\'efinie
ci-dessus en \'etant un exemple. En particulier, Panyushev donne un proc\'ed\'e pour trouver un syst\`eme de
g\'en\'erateurs de l'alg\`ebre
$Y(\mathfrak{k}')$, $\mathfrak{k}'$ \'etant une $\mathbb{Z}_2$-contract\'ee, et ce proc\'ed\'e est visible dans les
formules (1;6) de \cite{Ra-3} r\'e\'ecrites ci-dessus pour l'alg\`ebre de Lie $\mathfrak{b}$. Toutefois, il faut
signaler que $\mathfrak{b}$ ne co\"{i}ncide pas avec $\mathfrak{k}$. L'application $M : \mathfrak{b}
\longrightarrow
\mathfrak{k}$ d\'efinie par :
$$
        M(x,u,v^*) = \begin{pmatrix} 
x        &u\\
v^*        &0\\
\end{pmatrix}
$$
est un monomorphisme d'alg\`ebres de Lie dont l'image est un id\'eal de codimension 1 dans $\mathfrak{k}$.

\vskip 7mm
\noindent
\textbf{2.10.}~\underline{Remarque} :  

\smallskip
        Suivant la d\'emonstration de 2.4, on voit qu'on a le r\'esultat suivant : Toute fonction num\'erique continue sur
$\mathfrak{b}^*$, qui est $Ad^*_B(A)$-invariante, est automatiquement $Ad^*(B)$-invariante.

\smallskip
        Concernant les fonctions g\'en\'eralis\'ees, on peut montrer qu'une fonction g\'en\'eralis\'ee sur l'ouvert
$\widetilde{\Omega}$, qui est localement $A$-invariante, est automatiquement localement $B$-invariante.

\vskip 10mm
\begin{large}
\noindent
\textbf{3~Les cas $IO(n)$ et $ISO(n)$}
\end{large}

\bigskip 
\noindent
\textbf{3.1.}~On consid\`ere ici les groupes $C = IO(n)$ et $C^o = ISO(n)$, pr\'ecis\'ement : $C = O(n)
\underset{\rho}{\times} \mathbb{C}^n$, $C^o = SO(n)\underset{\rho}{\times} \mathbb{C}^n$, (avec $\rho(g).u = g.u$
dans les deux cas). L'alg\`ebre de Lie commune
\`a ces deux groupes est $\mathfrak{c} = so(n) \underset{\rho '}{\times} \mathbb{C}^n$. Les invariants dans
$S(\mathfrak{c})$ sous l'action coadjointe de $C$ et de $C^o$, not\'es respectivement $J$ et $Y(\mathfrak{c})$,
sont a priori distincts, avec $J \subset Y(\mathfrak{c})$. On montrera plus bas que $J = Y(\mathfrak{c})$
lorsque $n$ est pair, et que par contre lorsque $n$ est impair, $Y(\mathfrak{c})$ est un $J$-module de rang 2.

\vskip 7mm
\noindent
\textbf{3.2.}~On peut identifier $\mathfrak{c}$ \`a une sous-alg\`ebre de Lie de l'alg\`ebre de Lie
$\mathfrak{b}$ du num\'ero pr\'ec\'edent, au moyen de l'application $\gamma : \mathfrak{c} \longrightarrow
\mathfrak{b}$, d\'efinie par : $\gamma(x,u) = (x,u,-  {}^tu)$ (ici et ailleurs ${}^tX$ d\'esigne la transpos\'ee de
la matrice $X$). Le dual $\mathfrak{c}^*$ s'identifie alors au sous-espace de $\mathfrak{b}^*$ constitu\'e par les
formes lin\'eaires $\ell = (y,w^*,\xi)$, avec $y$ dans $so(n)$, et $\xi = - {}^tw^*$, de sorte qu'on peut parler de la
restriction \`a $c^*$ des fonctions d\'efinies dans $\mathfrak{b}^*$. On a en particulier : 
$$
        F_k(y, w^*, - {}^tw^*) = -w^*B_k(y) {}^tw^*
$$
Lorsque $k$ est impair, la matrice $B_k(y)$ est antisym\'etrique, et par suite, la restriction de $F_k$ \`a $c^*$
est nulle. On est amen\'e \`a consid\'erer les fonctions $\psi_k : c^* \longrightarrow \mathbb{C}$ d\'efinies par :
$$
        \psi_k(y,w^*) = -w^*B_{2k}(y){}^tw^*\qquad (0 \leq k \leq [\frac{n}{2}])
$$
(ici et ailleurs, $[\frac{n}{2}]$ est la partie enti\`ere de $\frac{n}{2}$).

\vskip 7mm
\noindent
\textbf{3.3.}~\underline{Remarques}
\vskip 1mm
\begin{enumerate}
\item Soit $\theta : \mathfrak{b} \longrightarrow \mathfrak{b}$ l'application lin\'eaire d\'efinie par :

$$
        \theta(x,u,v^*) = -({}^tx, {}^tv^*, {}^tu) \qquad (x,u,v^*) \in \mathfrak{b}
$$

qui n'est autre que l'oppos\'ee de la transposition des matrices lorsqu'on identifie les $(x,u,v^*)$ dans
$\mathfrak{b}$ aux matrices 
$\begin{pmatrix}
x        &u\\
v^*        &0\\
\end{pmatrix}$
(de taille $(n+1)$). On v\'erifie que $\theta$ est un automorphisme d'ordre 2 de l'alg\`ebre de Lie
$\mathfrak{b}$, et il est visible que $\gamma(\mathfrak{c})$ est l'espace des points fixes de $\theta$.  

\medskip
\item Soit $\cal{N}$ le sous-groupe de $Ad^*(B)$ qui normalise $\mathfrak{c}^* \subset \mathfrak{b}^*$. On
v\'erifie que $\cal{N}$ est l'ensemble des $Ad^*_B(g,u,- {}^tu)$, avec $g$ dans $O(n)$ et $u$ dans $V$, et que
l'image de $\cal{N}$ dans $GL(\mathfrak{c}^*)$ qui en r\'esulte co\"{i}ncide avec $Ad^*(C)$.
\end{enumerate}

\vskip 7mm
\noindent
\textbf{3.4.}~\textbf{Lemme} : \textit{Les fonctions $\Psi_k$ sont $Ad^*(C)$-invariantes.}

\vskip 2mm

        L'invariance des fonctions $\Psi_k$ r\'esulte imm\'ediatement de celle des fonctions $F_k$ sous l'action
coadjointe de $B$, et de la remarque (2) ci-dessus.

\vskip 7mm
\noindent
\textbf{3.5.}~\underline{Note bibliographique} :

\vskip 1mm
        Dans \cite{Ka-2}, Kaneta construit un sous-espace vectoriel $\mathfrak{h}$ de $\mathfrak{c}^*$, essentiellement
:
$\mathfrak{h} = \mathfrak{h}_0 \times \mathbb{C}\, e^*_n$, o\`u $\mathfrak{h}_0$ est un tore maximal de
$so(n-1)$ ($so(n-1)$ \'etant plong\'e naturellement, ``en haut et \`a gauche'', dans $so(n)$). Explicitement
(apr\`es complexification)
$$
\mathfrak{h}_0 = \{z = 
\begin{pmatrix}
0        &a_1        &        &        &        &        0 \\
-a_1        &0        &        &        &        &        0 \\
        &                &        \ddots        &        &        &        \vdots \\
        &        &        &        0        &a_\ell        &0 \\
        &        &        &        -a_\ell        &0        &0 \\
0        &0        &\ldots        &0        &0        &        0 \\
\end{pmatrix} ; \quad (a_1,a_2,\ldots ,a_\ell) \in \mathbb{C}^\ell \}, \hbox{lorsque}\ n=2\ell +1
$$

$$
\mathfrak{h}_0 = \{z = 
\begin{pmatrix}
0        &a_1        &        &        &        &        0&        0 \\
-a_1        &0        &        &        &        &        0&        0 \\
        &                &        \ddots        &        &        &        \vdots        &        \vdots \\
        &        &        &        0        &a_\ell        &0        &0 \\
        &        &        &        -a_\ell        &0        &0         &0\\
0        &0        &\ldots        &0        &0        &        0         &0\\
0        &0        &\ldots        &0        &0        &        0         &0\\
\end{pmatrix} ; \quad (a_1,a_2,\ldots ,a_\ell) \in \mathbb{C}^\ell \}, \hbox{lorsque}\ n=2\ell +2
$$
Il est facile de voir que l'op\'eration de restriction \`a $\mathfrak{h}$ des fonctions polyn\^omes appartenant
\`a $Y(\mathfrak{c})$ est un monomorphisme d'alg\`ebres dans l'alg\`ebre des fonctions polyn\^omes sur
$\mathfrak{h}$ (\textsl{c.f.} Lemma 3.5 dans \cite{Ka-2}, o\`u le corps de base est $\mathbb{R}$). Ceci \'etant,
Kaneta d\'efinit des fonctions $\varphi_0, \varphi_1,\ldots , \varphi_\ell$ sur $\mathfrak{h}$, de la mani\`ere
suivante :

$$
\begin{array}{ll}
\varphi_k(z,ae^*_n) = a^2\sigma_k(a^2_1,\ldots , a^2_\ell)         &(0 \leq k \leq \ell - 1) \\
\varphi_\ell(z,ae^*_n) = a a_1 a_2 \ldots a_\ell                &\hbox{lorsque} \ n = 2\ell +1 \\
\varphi_\ell(z,ae^*_n) = a^2\sigma_\ell(a^2_1,\ldots , a^2_\ell)                  &\hbox{lorsque} \ n = 2\ell +2 \\
\end{array}
$$

o\`u $\sigma_k\ (0\leq k \leq \ell)$ est la $k^{\hbox{i\`eme}}$ fonction sym\'etrique \'el\'ementaire
$(\sigma_0 = 1)$ ;

\vskip 2mm
et d\'emontre le : 

\vfill\eject 
%\vskip 7mm
\textbf{Th\'eor\`eme} (Theorem 2 dans \cite{Ka-2})

\vskip 1mm
        \textit{La $\mathbb{C}$-alg\`ebre $Y(\mathfrak{c})$ est isomorphe, via l'application de restriction $\psi
\longmapsto \psi|_\mathfrak{h}$, \`a la $\mathbb{C}$-alg\`ebre $\mathbb{C}[\varphi_0,\ldots , \varphi_\ell]$. Les
polyn\^omes
$\varphi_0, \varphi_1, \ldots , \varphi_\ell$, sont alg\'ebriquement ind\'ependants. }

\vskip 7mm
\noindent
\textbf{3.6.}~\textbf{Proposition} : \textit{Lorsque $n=2\ell + 2$, l'alg\`ebre $J$ des invariants de $IO(n)$
dans $S(\mathfrak{c})$ co\"{i}ncide avec $Y(\mathfrak{c})$ et :
$$
        J = \mathbb{C}[\psi_0,\ldots , \psi_\ell]
$$}

\noindent
\textsc{D\'emonstration} : Lorsque $(z,a\, e^*_n)$ appartient \`a $\mathfrak{h}$, il vient :
$$
        \psi_k(z,a\, e^*_n) = a^2 \sigma_k(a^2_1, \ldots , a^2_\ell) = \varphi_k(z,a\, e^*_n)
$$
Comme les $\psi_k$ sont $Ad^*(IO(n))$-invariantes, on voit que $Y(\mathfrak{c}) \subset J$. Par suite, les
invariants de $ISO(n)$ co\"{i}ncident avec les invariants de $IO(n)$.

\vskip 7mm
\noindent
\textbf{3.7.}~Supposons $n = 2\ell +1$. Comme vu plus haut, les fonctions $\varphi_0, \varphi_1,\ldots ,
\varphi_{\ell -1}$, sont les restrictions \`a $\mathfrak{h}$ des fonctions $IO(n)$-invariantes $\psi_0,\ldots ,
\psi_{\ell -1}$. Par contre, $\varphi^2_\ell$ est la restriction \`a $\mathfrak{h}$ de la fonction $\psi_\ell$. 

\vskip 1mm
        Soit $g = \begin{pmatrix} I_{2\ell}        &0 \\        0        &-1 \\                \end{pmatrix}$ (o\`u $I_{2\ell}$ est la matrice unit\'e de
taille $2\ell$). On a : $Ad^*(g,0).(z,ae^*_n) = (z,-ae^*_n)$, de sorte que : 
$$
\begin{array}{ll}
        \varphi_k(Ad^*(g,0).(z,ae^*_n)) = \varphi_k(z, ae^*_n) &(0 \leq k \leq \ell -1) \\
\varphi_\ell(Ad^*(g,0).(z,ae^*_n)) = -\varphi_\ell(z, ae^*_n) & \\
\end{array}
$$
Par suite : 
$$
        J = \mathbb{C}[\psi_0,\ldots , \psi_\ell]
$$
$$
        Y(\mathfrak{c}) = \mathbb{C}[\psi_0,\ldots , \psi_{\ell-1}, \Phi]
$$        
o\`u $\Phi$ est l'unique polyn\^ome dans $Y(\mathfrak{c})$ tel que : $\Phi|_\mathfrak{h} = \varphi_\ell$.

\vskip 7mm
\noindent
\textbf{3.8.}~Il reste \`a expliciter le polyn\^ome $\Phi$. Comme : 
$$
        \Phi^2(y,w^*) = \psi_\ell(y,w^*) = -w^*B_{2\ell}(y){}^tw^*
$$
on voit que $\Phi^2(y,w^*) = det\, Y$, avec $Y = \begin{pmatrix} y        &-{}^tw^* \\        w^*        &0 \\                \end{pmatrix}$.
Donc $\Phi$ co\"{i}ncide, au signe pr\`es, avec le pfaffien $Pf(Y)$ de la matrice $Y$ (de taille $(2\ell +2)$).

\vskip 1mm
Ici, il est peut-\^etre utile de faire le lien avec ce qui est appel\'e ``le pfaffien vectoriel'' $pf : so(2\ell +1)
\longrightarrow \mathbb{C}^{2\ell +1}$ dans (\cite{Ra-3}, \no 2.6). Il s'agit d'un covariant polynomial qui est li\'e au
pfaffien scalaire, de la mani\`ere suivante :
$$
        w^*pf(y) = Pf\begin{pmatrix} y        &-{}^tw^* \\        w^*        &0 \\                \end{pmatrix}
$$
et : $pf(gyg^{-1}) = (det\, g)g.pf(y)\quad (g \in\, O(2\ell +1))$.

\vskip 1mm
        Ainsi, l'invariant ``exotique'' $\Phi$ est, au signe pr\`es : 
$$
        \Phi(y,w^*) = w^*pf(y).
$$

\vskip 7mm
\noindent
\textbf{3.9.}~\underline{Note bibliographique} :

\vskip 1mm
        Dans \cite{Pa}, les invariants polyn\^omes pour les $\mathbb{Z}_2$-r\'eductions not\'ees 
$$
(\mathfrak{g}, \mathfrak{g_0}) = (so(n+m), so(m) \oplus so(n))
$$
sont d\'etermin\'es (Theorem 4.4) par une m\'ethode g\'en\'erale. Le cas $ISO(n)$ \'etudi\'e ci-dessus en est un
exemple.

%\vfill\eject
\vskip 7mm
\noindent
\textbf{3.10.}~\textbf{Remarques} : 
\vskip 1mm
\begin{enumerate}
\item 
Comme dans 2.8, il est possible d'avoir d'autres expressions pour les polyn\^omes $\psi_k$, du type : 
$$
        \psi_k(y,w^*) = p_{2k+2}(Y) - p_{2k+2}(y)
$$
avec : $Y = \begin{pmatrix} y        &-{}^tw^* \\        w^*        &0 \\                \end{pmatrix}$, et il est possible d'en d\'eduire que le
commutant de $so(n)$ dans l'alg\`ebre enveloppante de $\mathfrak{c}$ est commutatif (\textsl{c.f.} \cite{Ra-3},
formules (1;7) et $\no 2.6$). 

\medskip
\item
Compte-tenu de ce qui pr\'ec\`ede, l'alg\`ebre $J$ des invariants de $IO(n)$ est l'image de $Y(\mathfrak{b})$
par l'op\'eration de restriction \`a $\mathfrak{c}^*$ des fonctions polynomiales sur $\mathfrak{b}^*$. Ainsi
$J$ est un quotient de $Y(\mathfrak{b})$. Il est probable qu'il y ait une explication ``conceptuelle'' \`a ce fait.

\end{enumerate}

%%%%%%%%%%%%%%%%%%%%%%%
\vskip 15mm

%\vfill\eject
\renewcommand{\refname}{Bibliographie}     %%%pour m\'emoire : \cite{ } dans le texte

\medskip

Mustapha RAIS

E-mail : mustapha.rais@math.univ-poitiers.fr

\end {document}